\documentclass[letterpaper, 10 pt, conference]{ieeeconf}  

\IEEEoverridecommandlockouts                              

\usepackage{xcolor}
\usepackage{soul}
\usepackage{amsmath}
\usepackage{amssymb}
\usepackage{physics}
\usepackage{siunitx}
\usepackage[utf8]{inputenc}
\usepackage{fancyhdr}
\usepackage{graphicx}
\usepackage[export]{adjustbox}
\usepackage[T1]{fontenc}
\usepackage{tabularx,ragged2e,booktabs}
\usepackage{caption}
\usepackage{subcaption}
\usepackage{float}
\usepackage{multirow}
\usepackage{tabularx,booktabs}
\usepackage{verbatim}
\usepackage{cite}
\usepackage{enumerate}
\usepackage{tabu}
\usepackage[utf8]{inputenc}
\usepackage{dtk-logos} 
\newcolumntype{C}{>{\centering\arraybackslash}X} 
\usepackage{lipsum}
\newcolumntype{L}{>{\RaggedRight\arraybackslash%
    \hangafter=1\hangindent=1em}X} 

\title{\LARGE \bf
A Constraint Handling Approach with Guaranteed Feasibility for Surrogate Based Optimization
}

\author{Ahmed Abouhussein$^{1,4}$, Nusrat Islam$^{2,4}$ and Yulia T. Peet$^{3,4}$
\thanks{*This work was supported by NSF CMMI Award No. 1762827}
\thanks{$^{1}$Ph.D. student, email:
        {\tt\small aabouhus@asu.edu}}%
\thanks{$^{2}$Research Assistant, email: 
        {\tt\small nislam4@asu.edu}}%
\thanks{$^{3}$Assistant Professor, email:
        {\tt\small ypeet@asu.edu}}%
\thanks{$^{4}$School of Engineering, Matter, Transport and Energy, Arizona State University}%
}

\newcommand{\RomanNumeralCaps}[1]
    {\MakeUppercase{\romannumeral #1}}

\begin{document}

\maketitle
\thispagestyle{empty}
\pagestyle{empty}


\begin{abstract}
Gradient-free optimization methods, such as surrogate based optimization (SBO) methods, and genetic (GAs), or evolutionary (EAs) algorithms
have gained popularity in the field of constrained optimization of expensive black-box functions. However, constraint-handling methods, by both classes of solvers, do not usually guarantee strictly feasible candidates during optimization. This can become an issue in applied engineering problems where design variables must remain feasible for simulations to not fail. We propose a constraint-handling method for computationally inexpensive constraint functions which guarantees strictly feasible candidates when using a surrogate-based optimizer. We compare our method to other SBO, GA/EA and gradient-based algorithms on two (relatively simple and relatively hard) analytical test functions, and an applied fully-resolved Computational Fluid Dynamics (CFD) problem concerned with optimization of an undulatory swimming of a fish-like body, and show that the proposed algorithm shows favorable results while guaranteeing feasible candidates. 
\end{abstract}


\section{INTRODUCTION}

Global optimization  accounts for the majority of problems in practical and engineering optimization.
In terms of global optimization, gradient-based approaches can have faster convergence rates when compared to gradient-free approaches. However, they are more likely to remain trapped at a local minimum. Additionally, a wide variety of gradient-based algorithms become inefficient when dealing with many real-world applications in the fields of science, medicine, and engineering design. The inefficiency arises when computing the objective function derivative becomes expensive or sometimes even in feasible \cite{c0}. Derivative-free approaches such as evolutionary algorithms (EAs) or genetic algorithms (GAs) offer a lucrative alternative considering their ability in not requiring any assumption on the objective function landscape or its derivatives. However, a high number of function evaluations required by such algorithms can be prohibitive. Attempts at alleviating the mentioned concerns gave rise to the field of a surrogate-based optimization \cite{c1}. This branch of optimization techniques makes use of a surrogate model built with the help of true function evaluations, when finding local or global optima. Surrogate management framework (SMF) is an example of such a technique which has been successfully applied to various engineering design problems \cite{c2,c3}.

The goal of this study is to propose a surrogate based optimization (SBO) algorithm capable of ensuring strictly feasible candidates during optimization while keeping the number of function evaluations low. The purpose behind having such an algorithm would be to solve practical engineering problems such as finding the optimum modes of locomotion for a soft-robot fish for efficient underwater propulsion. 

There exists a wide spectrum of algorithms under the umbrella of SBO with varying performance depending on choices made in developing the surrogate optimization framework \cite{c4,c5}. The choices can vary in regards of the model used to construct the surrogate, the initial design of experiments sampling strategy and the infill strategy used to refine the surrogate over the optimization cycle. As the focus of this paper is feasible candidates, we restrict our discussion to the infill strategies as well as the constraint handling methods associated with SBO algorithms. 

One known infill strategy, proposed by Jones et al. \cite{c6}, provides a good balance between global and local search by choosing a candidate which maximizes expected improvement of the surrogate model. Another technique relies on a measure of surrogate "bumpiness" as it searches for a candidate \cite{c7}. The idea, in rough terms, is to provide an estimate of the objective function minimum and choose a candidate which is likely to produce that target objective while creating the least "bumpy" surrogate. A third technique, proposed by Regis and Shoemaker \cite{c8}, involves creating two groups of candidate points. The first group is generated by adding normally disturbed perturbations to the best found minimizer in a small vicinity close to the minimizer. The second group is generated by uniformly generating points from a box constrained domain of solutions. This technique will be explained in thorough detail later on as it forms the basis of the proposed constraint handling method.

As with infill strategy, there exists a variety of constraint handling methods in SBO algorithms. For example, a popular class of methods involves adding a penalty term to the objective function \cite{c9}. The penalty term seeks to assign a high cost for constraint violations, hence driving the algorithm away from infeasible candidates. A different approach considered for box-constrained problems is to recast infeasible candidates back into the feasible regime by reflecting candidates across the constraint boundary \cite{c10}. Lastly a simple rejection mechanism can be used to reject infeasible candidate. We propose an example of such a method and show that it can be applied to linear and nonlinear computationally inexpensive constraints. 

Our surrogate-based algorithm will be compared alongside various gradient-based and gradient-free algorithms (Table 1) on a three test problems: 1) the analytical Rosenbrock function \cite{c11}, 2) the Shifted Rotated Rastrigin’s Function \cite{c12}, 3) optimization of the locomotion of a thunniform bio-inspired propulsor based on Computational Fluid Dynamics (CFD) results. Algorithm performance metrics include a function evaluation count and a solution error to assess accuracy and convergence speed.  A formal problem presentation is given next in Section 2. A method description of the SBO algorithm as well as the constraint handling approach is presented in Section 3. Details regarding the computational set-up of the problems and the results are given in Section 4. Finally a  discussion of the results follows in Section 5.

\begin{table*}[h]
\centering
\caption{Optimization Algorithms}
\begin{tabular}{c c c c}
\toprule
 & Software Package & Solver & Comments \\ 
\midrule
\multirow{5}{*}{Gradient-based}     & MATLAB        & Interior Point (IP)                                                             & FMINCON\\ 
                                    & MATLAB        & Sequential Quadratic Programming (SQP)                                          & FMINCON\\ 
                                    & MATLAB        & Active Set (AS)                                                                 & FMINCON\\ 
                                    & DAKOTA        & Method of Feasible Directions (MFD)                                             & CONMIN\\ 
\hline
\multirow{9}{*}{Gradient-free}      & MATLAB        & Genetic Algorithm (GA)                                                          & \\ 
                                    & MATLAB        & Covariance Matrix Adaptation Evolution Strategy (CMAE-ES)                       & Hansen \cite{c13}\\ 
                                    & DAKOTA        & Evolutionary Algorithm (EA)                                                     & COLINY \cite{c14}\\ 
                                    & MATLAB        & Surrogate Based Optimizer (SBO)                                                 & w/ constraint handling\\ 
                                    & MATLAB        & Surrogate Based Optimizer Penalty (SBO\_P)                                      & w/o constraint handling\\  
                                    & PYTHON        & DYnamic COordinate search using Response Surfaces (DYCORS) from pySOT           & DYCORS \cite{c10}, pySOT \cite{c15}: \\
                                    & DAKOTA        & Surrogate based global optimizer (SBG)                                          & \\ 
                                    
\bottomrule
\end{tabular}
\end{table*}


\section{PROBLEM FORMULATION}
The general formulation of the optimization problem considered can be expressed as follows,
\begin{equation} 
\begin{split}
  & \text{minimize} \quad f(\mathbf{x}) \linebreak \\
  & \text{subject to} \quad \mathbf{x} \in \mathbb{R}^n
  \end{split}
\end{equation}
where $f:\mathbb{R}^n \rightarrow \mathbb{R}$ is the objective function, $\mathbf{x}$ is a vector of design parameters. The set $S\subseteq \mathbb{R}^n$ contains the $n$-dimensional search space, which would define a rectangle in $\mathbb{R}^2$ or a rectangular cuboid in $\mathbb{R}^3$:
\begin{equation} 
l(i) \leq x_{i} \leq u(i), \quad 1 \leq i \leq n
\end{equation}
where $l(i)$ and $u(i)$ represent lower and upper bounds, respectively, on a design parameter in the $i$th dimension. The set $C\subseteq \mathbb{R}^n$ contains a set of $m \geq 0$ constraints:
\begin{equation} 
  \begin{split}
  & g_{r}(\mathbf{x}) \leq 0, \quad r = 1, ..., q,\\
  & h_{r}(\mathbf{x}) = 0, \quad r = q + 1, ..., m
  \end{split}
\end{equation}
where $g_{r}(\mathbf{x})$ and $ h_{r}(\mathbf{x})$ are referred to as the inequality and equality constraint sets, respectively, on the design parameter vector, $\mathbf{x}$.
\subsection{Test Problem \RomanNumeralCaps{1}}
The Rosenbrock function is a canonical optimization test function known for containing a global minimum within a wide basin. The two-dimensional form of the function used in this study is given as follows:
\begin{equation}
f(x_{1},x_{2})= (a-x_{1})^{2}+b(x_{2}-x_{1}^{2})^{2}
\end{equation}
with global minimum located at $(x^*_{1},x^*_{2})=(a,a^2)$, where $f(x^*_{1},x^*_{2})=0$. If $\mathbf{x} \not\in S \cap C$, the function is given a penalty value of $\num{1e+6}$. The parameters were chosen to be $a=0.35$ and $b=100$ to insure that the global minimum resides within the constraint domain as shown in Fig. 1.  The $C$ constraint set, chosen to resemble the natural constraint in the third test problem, is given below:
\begin{equation} 
  \begin{split}
  & g_{1}(\mathbf{x}) = x_{2}+2.5x_{1}^2-0.5 \leq 0 \\
  & g_{2}(\mathbf{x}) = -x_{2}-x_{1} + 0.4 \leq 0
  \end{split}
\end{equation}
\begin{figure}[h]
    \centering
    \includegraphics[width=0.47\textwidth, inner]{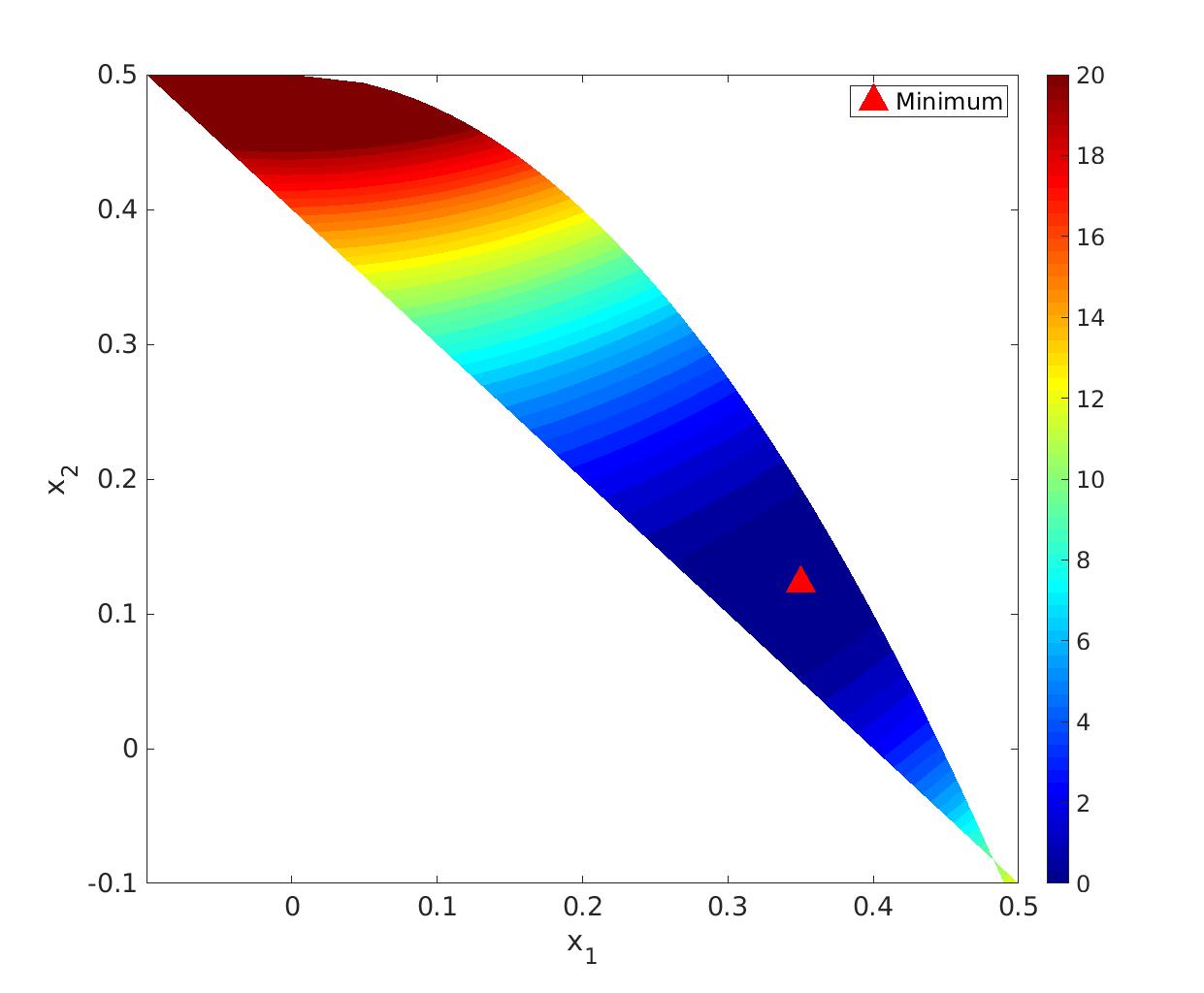}
    \caption{2-D Rosenbrock over the constrained domain}
\end{figure}
\subsection{Test Problem \RomanNumeralCaps{2}}
The second analytical function considered is the shifted rotated Rastrigin function used as one of the benchmark functions for CEC2005 competition\cite{c12}. This variant of the Rastrigin function provides a non-linear, non-separable, highly multi-modal challenging test function (Figure 2). The 2-D form used in this study is given as follows:
\begin{equation}
f(x_{1},x_{2})= \sum_{i=1}^{2}(z_{i}^{2}-10\cos(2\pi z_{i}))-330
\end{equation}
where $\mathbf{z} = (\mathbf{x}-\mathbf{o})\mathbf{M}$, $\mathbf{x} = [x_{1},x_{2}]$, $\mathbf{o} = [o_{1},o_{2}]$ is the shifted global optimum, $\mathbf{M}$ is a linear transformation matrix with condition number = 2. As with Test Problem 1, the function was given a penalty value of $\num{1e+6}$ in the unfeasible domain. The $C$ constraint set in this case is:
\begin{equation} 
  \begin{split}
  & g_{1}(\mathbf{x}) = x_{2}+2.5x_{1}^2-0.5 \leq 0 \\
  & g_{2}(\mathbf{x}) = -x_{2}-x_{1} + 0.4 \leq 0
  \end{split}
\end{equation}
\begin{figure}[h]
    \centering
    \includegraphics[width=0.47\textwidth, inner]{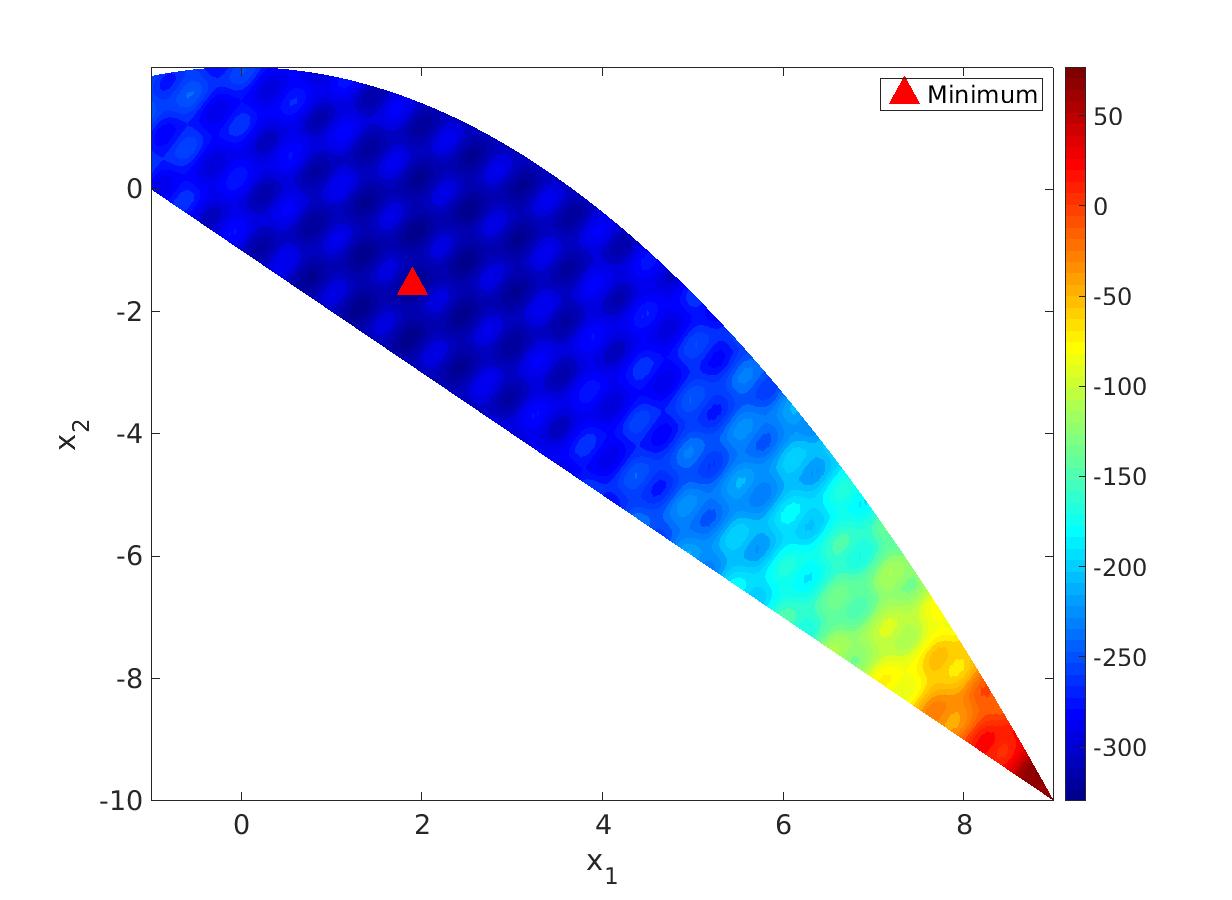}
    \caption{2-D Shifted Rotated Rastrigin over the constrained domain}
\end{figure}
\subsection{Test Problem \RomanNumeralCaps{3}}
The black-box simulation is based on the modeling and computational framework of a 2D single thunniform fish swimming under water presented by Xu and Peet \cite{c16}. The optimization process aims to locate an optimum fishing mode related to a kinematic gait that maximizes the start-up propulsive force with minimal energy expended, referred to as the propulsive efficiency. The fish undulation is described in terms of its center-line position:
\begin{equation} 
y_{c}(p,t) = [x_1 \frac{p}{L} + x_2 (\frac{p}{L})^2]\sin(2\pi(\frac{p}{\lambda L} - ft))
\end{equation}
where $x_{1}$ and $x_{2}$ are dimensional undetermined linear and quadratic wave amplitude coefficients, $p$ is the streamwise fish position, $y$ is the spanwise fish position, $\lambda$ is the body wave length which is 1.1, $f$ is the body wave tail-beat frequency taken to be 1 Hz, $L$ is the dimensional length of the fish taken to be $0.3$ m and $t$ is time.  Consequently, the coefficients $\{ x_{1}, x_{2} \}$ span a range of swimming modes.  In order to allow for only physically realizable modes, the following constraint set was imposed \cite{c17,c18,c19}:
\begin{equation} 
  \begin{split}
  & g_{1}(\mathbf{x}) = 0.4x_{2}L + x_{1}^2 \leq 0 \\
  & g_{2}(\mathbf{x}) = |x_{2}+x_{1}| - 0.1L \leq 0 \\
  & g_{3}(\mathbf{x}) = x_{2} \leq 0
  \end{split}
\end{equation}
The propulsive efficiency, $\eta$, is defined as follows: 
\begin{equation} 
\eta(x_{1},x_{2},t) = \frac{\int_{0}^{t} \oint_{body} -\sigma \cdot n_{p} \cdot U \dd p \dd t}{\int_{0}^{t} \oint_{body} -\sigma \cdot n_{y} \cdot \nu \dd p \dd t}
\end{equation}
where $\sigma$ is he total Cauchy stress tensor which includes viscous and pressure forces, $n = \{ n_{p},n_{y}\}$ is the outer unit normal vector on the body surface, $\nu(p,t) = {\partial y_{m}(p,t)}/{\partial t}$ is the surface transverse velocity due to undulation, and $U(t)$  is the propulsive forward velocity.  During optimization, the negative value of the propulsive efficiency is taken as the objective function and the $\{ x_{1}, x_{2} \}$ coefficients are taken as design parameters subject to the feasible parameter space.
\section{Method Description}
SBO relies on approximating solutions based on a surrogate model of the objective function. First an initial surrogate model of the objective function is created using a data set of true function evaluations sampled with a space filling strategy. Then, with each iteration an infill criteria, which attempts to balance global and local exploration, is used to refine the surrogate with multiple surrogate function evaluations and one true function evaluation. In this work, we choose a sampling technique and a surrogate model from the DACE \cite{c20} Matlab tool box. We define an infill strategy according to the Metric Stochastic Response Surface (\textit{MSRS}) method proposed by Regis and Shoemaker\cite{c8}. Finally, we extend the \textit{MSRS} method to produce strictly feasible candidates at each iteration step.

We present below the SBO used in this study. Let the feasible domain, $D$, be defined as $D = S \, \cap \, C$. Let $M_k$ and $s_k(x)$ be defined as the mesh space and the surrogate model, respectively, at iteration $k$. Define maximum number of iterations, $k_{max}$, and a tolerance, $tol$.
\begin{enumerate}[\textit{Step} 1:]
    \item Sample a finite set of evaluation points $T \subset D$ using Latin Hypercube Sampling\cite{c21}, where $\text{card}(T)=20$. Evaluate $f(T)$, where $f(x)$ is the true objective function. Identify the current best point $x_0$. Set $M_0 = T$.
    \item Fit a Kriging model surrogate, $s_0(x)$, with a Gaussian correlation function and $0$-th order regression polynomial. In other words, the surrogate model is assumed to have a constant mean and a stochastic error term that is modeled by a Gaussian process. This is referred to as ordinary kriging and allows for a flexible and reliable prediction method \cite{c22}.
    \item While $(k < k_{max})$
        \begin{enumerate}
        \item Create a set of strictly feasible candidate points, $X_k$, according to the proposed algorithm (section 3.A.) and evaluate $s_k(X_k)$. 
        \item Use the \textit{MSRS} method which assigns a weighted score to each point in set $X_k$ based on two criteria: 1) the distance of points in $X_k$ to $M_k$, and 2) the surrogate response values, $s_k(X_k)$. The weighted score insures that the next candidate point has a low objective value that is far away from already sampled points. The point with the best weighted score is identified as the next evaluation point, $x_k$. 
        \item Evaluate $f(x_k)$. If $tol$ is met: \textbf{break}.
        \item Set $M_{k+1} = M_k + x_k$. Re-fit $s_{k+1}$ with $M_{k+1}$. Set $k = k + 1$.
        \end{enumerate}
    \item Return $x_k$
\end{enumerate}
\subsection{Constrained Candidate Sampling}
The candidate points are split into two categories\cite{c8}:
\begin{enumerate}
    \item Uniformly sampled global points: The first set, $U_k$, is generated by a uniform random sampling of points from the box-constrained domain such that $U_k \subset S$. We set $\text{card}(U)=2000$.
    \item Normally sampled local points: The second set $N_k$ is generated by adding perturbations to $x_k$ drawn from a random normal distribution with zero mean and unit variance. The are three perturbation rates chosen: one-tenth, one-hundredth and one-thousandth of the smallest variable range. The smallest variable range is defined as $\text{min}(u-l)$. We set $\text{card}(N)=2000$.
\end{enumerate}
We define the possibly unfeasible candidate set as $X'_k = U_k + N_k$. It follows that $\text{card}(X'_k)=4000$. We note that $X'_k \subset S$, however $X'_k \not\subset S \cap C$. We then enforce all constraints, linear and/or nonlinear, through the following algorithm:
\begin{enumerate}[\textit{Step} 1:]
    \item Evaluate $g_{r}(X'_k)$ for $r=1, ..., q$.
    \item Define `penalty' vector, $J$, for each candidate point:
        \begin{equation} 
           j_{i} =\sum_{r=1}^{m}\text{max}(0,g_{r}(x_{i})), \quad i = 1, ..., \text{card}(X'_k)\\
        \end{equation}
        where $j_{i}$ are the entries of set $J$.
    \item The candidate points set, $X_k$, is simply defined as: $X_k = X'_k( J = 0 )$. In other words, the $X_k$ candidates are the candidates in $X'_k$ with a zero penalty.
\end{enumerate}
This method does not guarantee that the candidate set will always be the of expected $\text{card}(X'_k)$ but the final candidate set $X_k$ is always feasible for all possible candidate solutions. If it happens that $X_k \subset M_k$, then the candidates are resampled.
\section{Computational Set-up \& Results}
\subsection{Test Problems \RomanNumeralCaps{1} \& \RomanNumeralCaps{2} }
\begin{table*}[h]
\caption{Rosenbrock Function Optimization Results}
\begin{tabularx}{\textwidth}{@{}l*{3}{C}c@{}}
\toprule
Algorithm     & Function Evaluation & Function Evaluation MOE & Absolute Error & Absolute Error MOE \\ 
\midrule

IP     & 24        & N/A    & $\num{6.95e-5}$      & N/A \\ 
SQP    & 26        & N/A    & $\num{3.59e-5}$      & N/A \\ 
AS     & 25        & N/A    & $\num{2.50e-5}$      & N/A \\ 
MFD    & 13        & N/A    & $\num{1.54e-6}$      & N/A \\ 
GA     & 1904      & 135    & $\num{6.43e-5}$      & $\num{2.04e-6}$ \\ 
CMAES  & 106       & 3      & $\num{4.97e-5}$      & $\num{1.85e-6}$ \\
EA     & 975       & 14     & $\num{7.49e-5}$      & $\num{5.53e-6}$ \\
SBO    & 32        & 0      & $\num{3.10e-5}$      & $\num{1.74e-6}$ \\
SBO\_P & 99        & 0      & $\num{8.24e+0}$      & $\num{4.73e-1}$ \\
DYCORS & 32        & N/A    & $\num{8.86e-2}$      & N/A \\
SBG    & 32        & N/A    & $\num{7.54e+0}$      & N/A \\

\bottomrule
\end{tabularx}
\end{table*}
\begin{table*}[h]
\caption{Shifted Rotated Rastrigin Function Optimization Results}
\begin{tabularx}{\textwidth}{@{}l*{3}{C}c@{}}
\toprule
Algorithm     & Function Evaluation & Function Evaluation MOE & Relative Error & Absolute Error MOE \\ 
\midrule

IP             & 68        & N/A    & $\num{3.92e-2}$        & N/A \\ 
SQP            & 7         & N/A    & $\num{1.46e-1}$        & N/A \\ 
AS             & 20        & N/A    & $\num{1.24e-1}$        & N/A \\ 
MFD            & 25        & N/A    & $\num{9.64e-2}$         & N/A \\ 
GA             & 1153      & 33     & $\num{2.12e-3}$         & $\num{4.93e-2}$ \\ 
CMAES          & 326       & 22     & $\num{5.06e-3}$         & $\num{8.91e-2}$\\
EA             & 253       & 12     & $\num{3.97e-3}$         & $\num{3.18e-2}$\\
SBO            & 226       & 12     & $\num{4.18e-3}$         & $\num{3.01e-2}$ \\
SBO\_P         & 466       & 25     & $\num{3.97e-3}$         & $\num{3.03e-2}$ \\
DYCORS         & 226       & N/A    & $\num{5.61e-3}$         & N/A \\
SBG            & 420       & N/A    & $\num{6.09e-2}$         & N/A \\

\bottomrule
\end{tabularx}
\end{table*}
Numerical simulations were performed on six 3.20GHz Intel processors on a Linux environment. The optimization algorithms were run on Matlab 2018b, Dakota 6.8 and Python 3.6.10. All solvers were started from the feasible point of $(x_{1},x_{2})_{0}=(0.2,0.3)$ and $(x_{1},x_{2})_{0}=(-0.3,0.5)$ and the convergence criteria was set to $\num{1E-4}$ and $\num{2E+0}$ on the absolute error whenever possible for Test Problems \RomanNumeralCaps{1} and \RomanNumeralCaps{2}, respectively. An exception to convergence termination criteria is the DYCORS and SBG solvers, which can only be terminated with a maximum iteration threshold. For these solvers we used the fastest function evaluation count from other algorithms as the maximum iteration threshold. We note that Test Problem \RomanNumeralCaps{2} is a highly challenging optimization problem and hence the error tolerance was adjusted accordingly. The function evaluations count as well as the error, achieved after convergence, are reported in Table 2 and Table 3. We use an absolute error metric for the first problem, where the true solution is 0, and a relative error metric for the second problem, where the true solution is -330. Additionally, a margin of error (MOE) is calculated with 95\% confidence for appropriate solvers.
as follows:
\begin{equation} 
  \begin{split}
   & SD = \sqrt{\frac{\sum_{i=1}^{m}(z_{i}-\mu_{z})^2}{N-1}} \\
   & MOE = \frac{SD}{\sqrt{N}}*1.96 \\
  \end{split}
\end{equation}
where $z_{i}$ is a random variable at iteration $i$, $\mu_z$ is the sample mean, $N$ is the number of realizations taken to be 1000, $SD$ is the standard deviation and $MOE$ is the margin of error. 
The exceptions are the deterministic gradient-based solvers, the SBG solver and DYCORS. The population size of all three EAs/GAs were chosen to be 20. The crossover rate and mutation rate were set to be 50\% for the COLINY\_EA solver, while those rates could not be specified for the other EA/GA algorithms. We note that penalty value of $\num{1e+6}$ would cause the CMAES solver to occasionally converge to that value with a small evaluation count. When averaging over $N$ realizations, the results showed CMAES to have a suspiciously low function evaluation count and an error on the order of $O(10^4)$. For this reason the penalty value assigned to this algorithm is the Matlab $NaN$, or “Not a Number”. The SBO algorithms used 20 function evaluations to construct an initial surrogate of the function. 
While, SBG used a derivative-free algorithm (COLINY\_EA) when minimizing, it was still unable to provide solutions within tolerance. In that case, an increased maximum iteration threshold had little effect on improving solution quality. 
\subsection{Test  Problem \RomanNumeralCaps{3}}
The simulations were performed on three nodes of a super computing cluster with each compute node containing two Intel Xeon E5-2680 v4 CPUs running at 2.40GHz. The COLINY\_EA solver and the MATLAB Surrogate Model Optimization \cite{c23} toolbox were chosen to solve the optimization problem. Simulations of the swimming fish were performed using Nek5000, a high-fidelity open source CFD solver\cite{c24}. The optimization was set to be terminated at a relative convergence tolerance of $\num{1e-6}$. The relative convergence tolerance is defined by the relative change of the objective function between successive iterations. 
In case of COLINY\_EA, the Dakota package managed the optimization process. We developed a batch script which accepts a set of test parameters from Dakota, updates and runs the Nek5000 CFD simulation accordingly, post-processes the simulation results and finally produces an output file which is read by Dakota. The SBO optimization followed a similar work flow however a MATLAB script was used to manage the optimization procedure. The final Kriging surrogate was used to create an objective function landscape contour found in Fig. 2. The results for this problem are shown in Table 4. 
\begin{figure}[ht!]
    \centering
    \includegraphics[width=0.47\textwidth, inner]{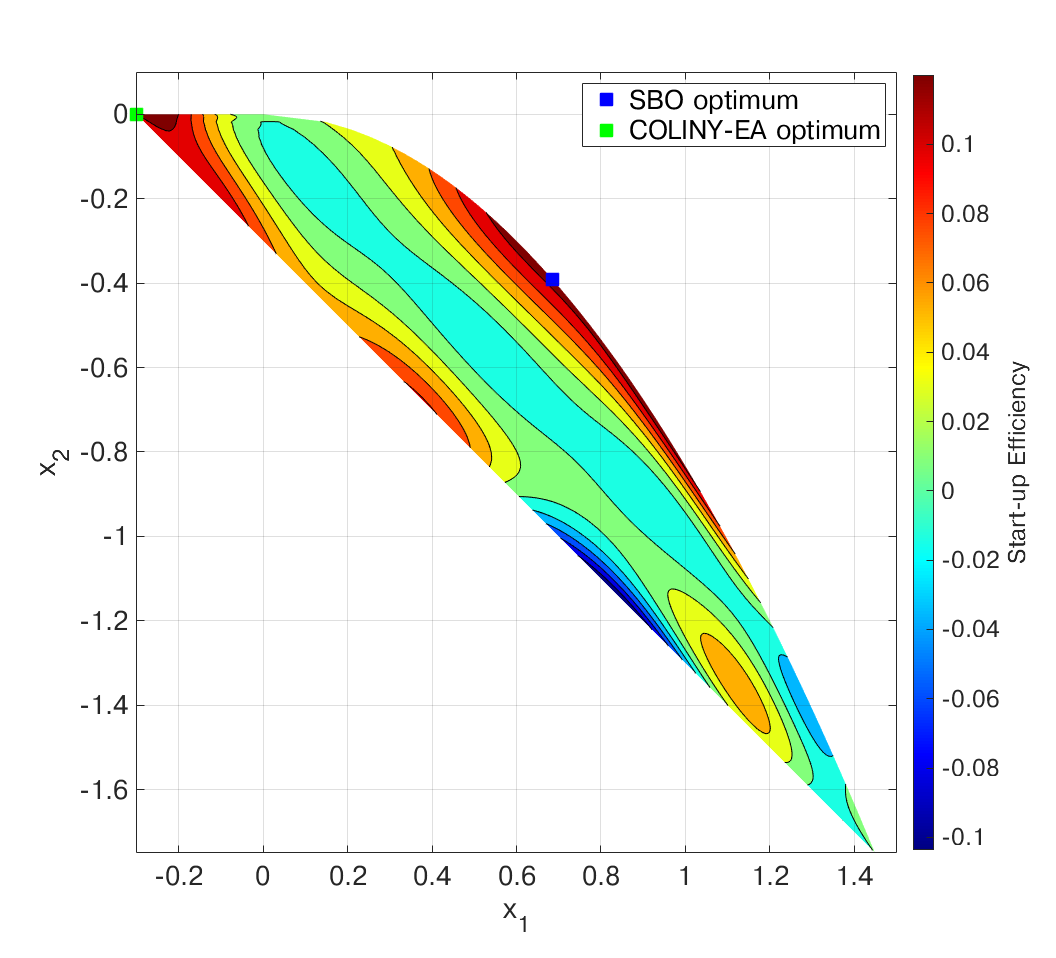}
    \caption{2-D start-up efficiency function over the constrained domain}
    \label{fig:galaxy}
\end{figure}
\begin{table}[h]
\caption{Black-box optimization}
\begin{tabularx}{0.5\textwidth}{@{}l*{4}{C}c@{}}
\toprule
Algorithm     & Function Evaluation & Start-up Efficiency & Optimized Parameter Set \{${x_{1},x_{2}}$\} \\ 
\midrule

COLINY\_EA     & 735       & $11.6\%$ & $\{-0.3 , 0\}$ \\
SBO            & 58        & $13.1\%$ & $\{0.68 ,-0.39\}$ \\

\bottomrule
\end{tabularx}
\end{table}

\begin{figure}[h!]  
    \begin{subfigure}{0.20\textwidth}
    \centering
    \includegraphics[width=\textwidth]{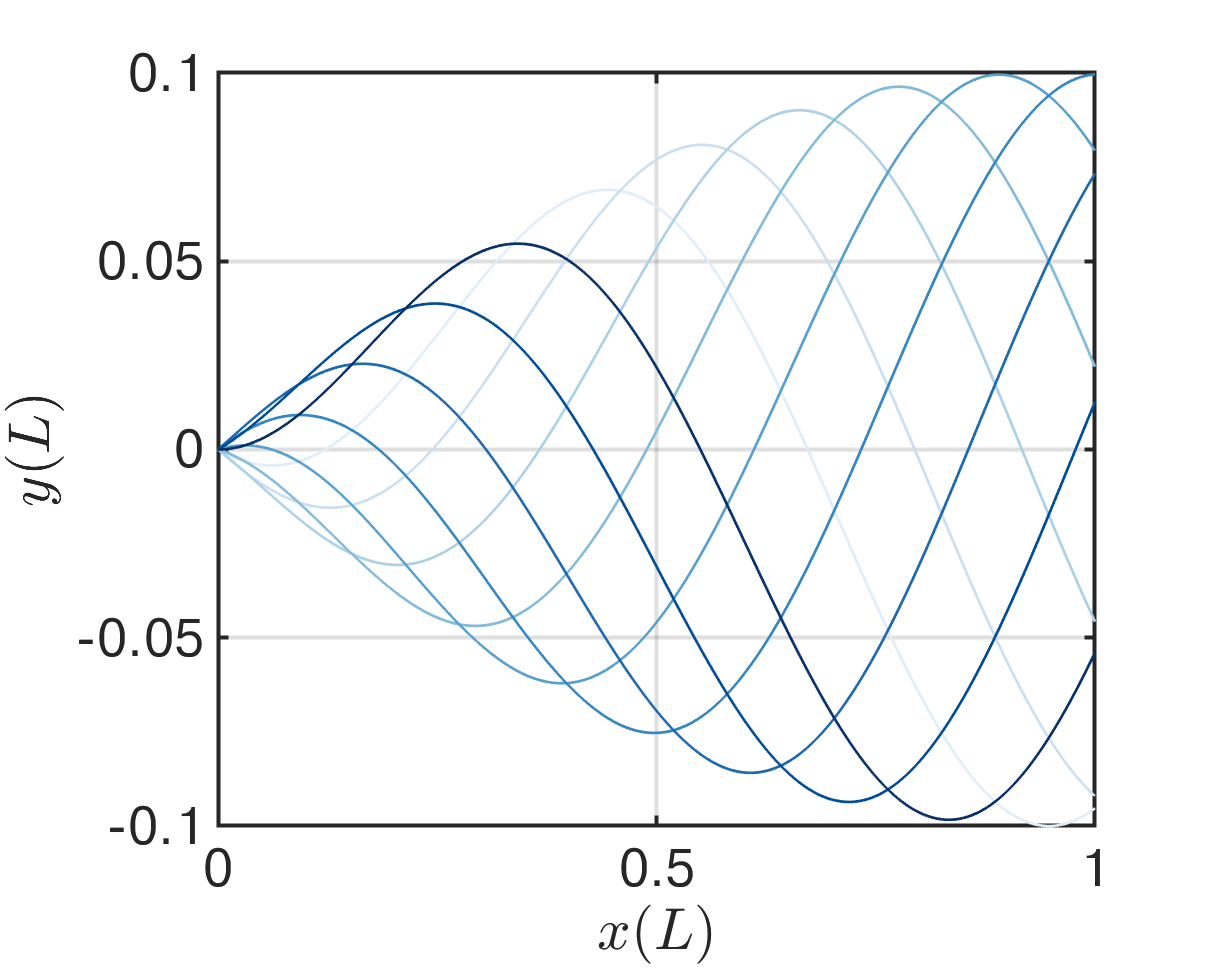}
    \caption{SBO Optimum}
    \end{subfigure}
     \begin{subfigure}{0.20\textwidth}
    \centering
    \includegraphics[width=\textwidth]{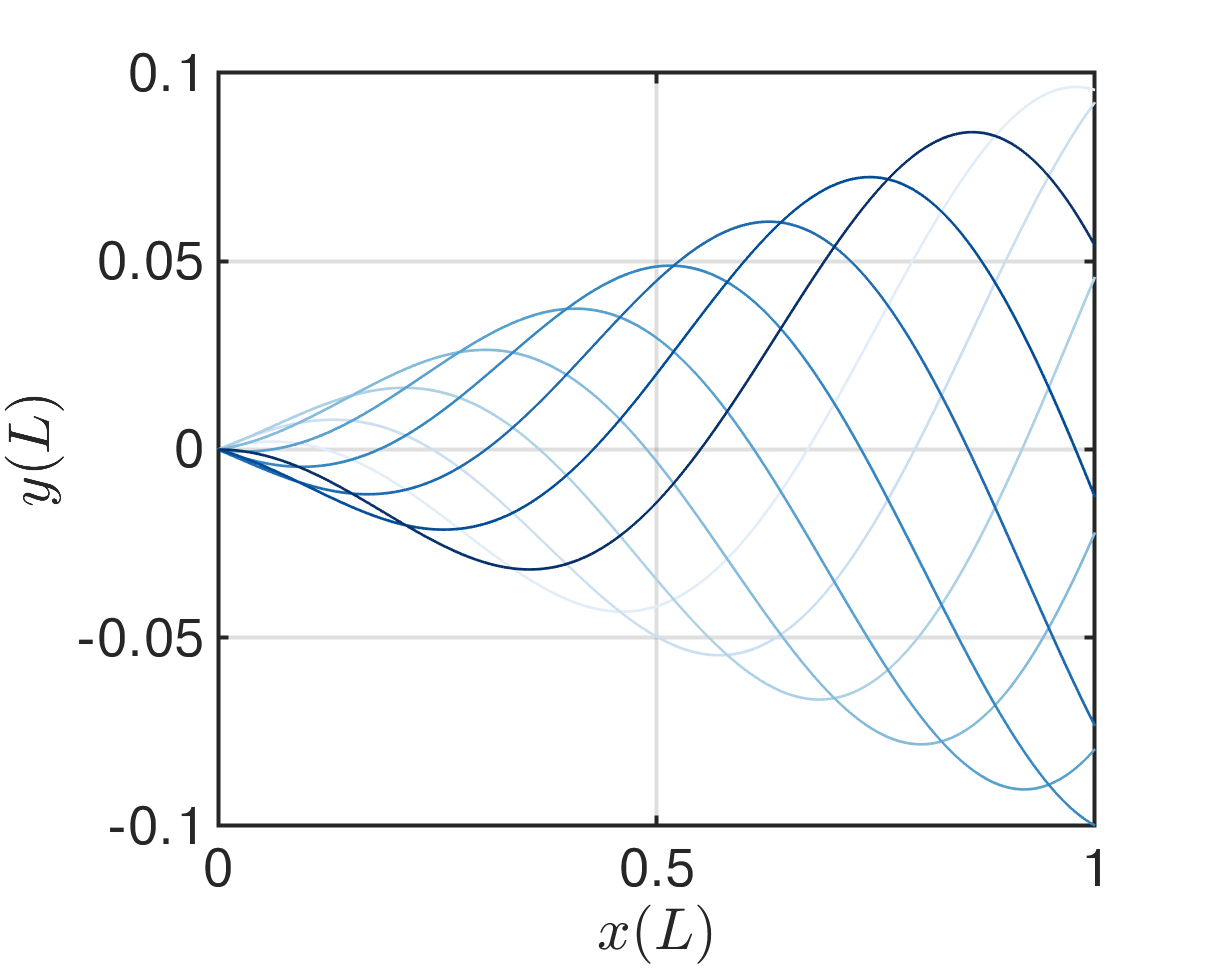}
    \caption{COLINY\_EA Optimum}
    \end{subfigure}
    \caption{Swimmer midline deformation across one time cycle for two propulsive modes: (a) SBO Optimum; and (b) COLINY\_EA Optimum. Deformations of the midline in time are encoded every $1/10^{\,th}$ of the period in different shades of blue from lightest ($t = 0$) to  darkest ($t = T$).}
    \label{fig:Fish position}
\end{figure}


\section{DISCUSSION}

The SBO algorithm with the proposed constraint handling approach has consistently outperformed other constraint handling algorithms on analytical test problems as well as an expensive black-box optimization problem. Simple box-constraint handling methods, such as DYCORS and SBO\_P, are shown to have less accuracy for a given number of function evaluations. While penalizing unfeasible candidates may not be effective with SBOs (such as SBG), it is shown to be an effective constraint-handling method in GAs (CMAES and GA). However, the first two problems showed that SBO can produce results on the same order of accuracy as GAs with as low as 30\% of the function evaluations needed for the latter. The third test problem, which considers an expensive black-box function, showed that the proposed SBO provides better global convergence, when compared with GAs, while only requiring as low as 8\% of the function evaluations needed for GAs. The global maximzer resulted in a notably different kinematic gait profile when compared to the local maximizer from the genetic algorithm (Fig. 3), resulting in changes in hydrodynamic quantities. We note that gradient-based methods, such as the algorithms available through FMINCON and CONMIN, should not be considered if the black-box optimization problem is not known to be unimodal. The results indicate that the proposed SBO constraint handling approach is a promising approach for dealing with constrained computationally expensive black-box problems which require strictly feasible candidates. 

\section{Acknowledgements}

This research is supported by NSF CMMI grant \# 1762827.


\addtolength{\textheight}{-12cm}   









\end{document}